\documentclass{amsart}
\usepackage{amssymb,mathrsfs,MnSymbol}
\usepackage{color}
\usepackage{tikz}
\usetikzlibrary{positioning}

\usepackage[textheight=25 cm, textwidth=16cm, headheight=10pt, headsep=20pt, footskip=10pt, left=1.5cm, right=1.5cm]{geometry}

\def\bR {\mathbf{R}}

\def\fH {\mathfrak{H}}

\def\cA {\mathcal{A}}
\def\cB {\mathcal{B}}
\def\cC {\mathcal{C}}
\def\cD {\mathcal{D}}

\def\\fH {\mathcal{H}}

\def\\hbar {{\\hbarilon}}

\newcommand{\Tr}{\operatorname{trace}}

\newcommand{\ba}{\begin{aligned}}
\newcommand{\ea}{\end{aligned}}

\newcommand{\be}{\begin{equation}}
\newcommand{\ee}{\end{equation}}

\newcommand{\bea}{\begin{eqnarray}}
\newcommand{\eea}{\end{eqnarray}}

\newcommand{\nn}{\nonumber}

\newcommand{\coup}{\Pi}

\newcommand{\ket}[1]{| #1\rangle}
\newcommand{\bra}[1]{\langle #1|}

\newcommand{\ketbra}[2]{\ket{#1}\bra{#2}}

\begin{document}

\title{\Large 
Quantum optimal transport is cheaper}

\author[\Large E. Caglioti]{Emanuele Caglioti}
\address[E.C.]{Sapienza Universit\`a di Roma, Dipartimento di Matematica, Piazzale Aldo Moro 5, 00185 Roma}
\email{caglioti@mat.uniroma1.it}

\author[F. Golse]{Fran\c cois Golse}
\address[F.G.]{CMLS, Ecole polytechnique,  91128 Palaiseau Cedex, France}
\email{francois.golse@polytechnique.edu}

\author[T. Paul]{Thierry Paul}
\address[T.P.]{CMLS, CNRS, Ecole polytechnique,   91128 Palaiseau Cedex, France}
\email{paul@ljll.math.upmc.fr}

\begin{abstract}
\Large
We compare bipartite (Euclidean) matching problems  in classical and quantum mechanics. The quantum case is defined after a quantum version of the Wasserstein distance introduced in \cite{GMouPaul}.  We show that the optimal quantum cost can be cheaper than the classical one. We treat in detail the case of two particles: the case of equal mass provides equal quantum and classical costs, while we exhibit examples of different masses for which the quantum cost is actually strictly cheaper. 
\end{abstract}
\LARGE
\maketitle
%
%
%
%
%
%
%

\tableofcontents
\section{Introduction}\label{intro}
The paradigm of modern optimal transport theory uses extensively the $2$-Wasserstein distance between two  probability measures $\mu,\nu$ on $\bR^n$, defined as
\be\label{claswass}
W_2(\mu,\nu)^2:=\inf_{\coup\text{ coupling of }\mu\text{ and }\nu}\int|x-y|^2\coup(dx,dy).
\ee
We have called  coupling
(or transport plan) of the two probabilities $\mu$ and $\nu$ any  probability measure $\coup(dx,dy)$ on $\bR^n\times\bR^n$ whose 
marginals on the first and the second factors are $\mu$ and $\nu$ resp., i.e. 
\be
\label{clamarg}
\begin{aligned}
&\int_{\bR^{n}\times\bR^n}a(x)\coup(dx,dy)=\int_{\bR^n}a(x)\mu(dx),\\
&\int_{\bR^n\times\bR^n}b(y)\coup(dx,dy)=\int_{\bR^n}b(y)\nu(dy)
\end{aligned}
\ee
 for all  test (i.e. continuous and bounded) functions $a$ and $b$.

\vskip 0.5cm
Restricting the definition of $W_2$ to couplings of the form 
\be\label{ksb}
\coup
=
\delta(y-T(x))\mu(dx)
\ee
 where $T$ is a  transformation of $\bR^n$ such that $\nu$ is the image
$T_{\#}\mu$ of $\mu$ by $T$, one sees that:
\be\label{mongeprob}
M(\mu,\nu)^2:=\inf_{T_{\#}\mu=\nu}\int_{\bR^n} (x-T(x))^2\mu(dx)\ge W_2(\mu,\nu)^2.
\ee
 
The converse inequality is due to Knott, Smith and Brenier: under certain restrictions on the regularity of $\mu$ and $\nu$, any optimal coupling for the minimization problem defined by \eqref{claswass} is of the form \eqref{ksb} 
for some transport map $T$, so that the inequality in \eqref{mongeprob} is an equality
 (see e.g. \cite{cgp} Section 1 for some details and 
Theorem 2.12 in \cite{VillaniTOT}  
for an extensive study).
\vskip 1cm
Associated to $W_2$ is the bipartite  matching problem which can be described as follows. 
Let us consider $M$ material points on the real line $\{x_i\}_{i=1,\dots,M}$ with $x_i<x_{i+1}$, and with masses $\{m_i\}_{i=1,\dots,M}$, and on the other hand $N$ points $\{y_i\}_{i=1,\dots,N}$ with $y _j<y_{j+1}$, and with masses 
$\{n_i\}_{i=1,\dots,N}$. We normalize the total mass as follows: 
$$
\sum\limits_{i=1}^Mm_i=\sum\limits_{j=1}^Nn_j=1.
$$
The bipartite problem consists in finding a coupling matrix $(p_{i,j})|_{i=1,\dots,N,j=1,\dots M}$ satisfying
$$
\sum_{j=1}^Np_{i,j}=m_i,\qquad \sum_{i=1}^Mp_{i,j}=n_j\,,\qquad p_{i,j}\ge 0\text{ for each }i,j
$$
which minimizes the quantity 
$
\sum\limits_{i,j}p_{i,j}|x_i-y_j|^2\,.
$

That is to say, we define the optimal transport cost as
$$
C
_c
:=\inf_{\substack{
p_{i,j}\ge 0\\ \sum\limits_{j=1}^Np_{i,j}=m_i,\,\,\,\sum\limits_{i=1}^Mp_{i,j}=n_j}} \sum\limits_{i,j}p_{i,j}|x_i-y_j|^2.
$$
It is natural to associate to the sets $\{x_i\}_{i=1,\dots,M}$ and $\{m_i\}_{i=1,\dots,M}$, and to the sets $\{y_i\}_{i=1,\dots,N}$ and $\{n_i\}_{i=1,\dots,N}$ the following discrete  probability measures
$$
\mu:=\sum_{i=1}^Mm_i\delta_{x_i},\qquad \nu:=\sum_{j=1}^Nn_j\delta_{y_j}.
$$
It is easy to see that any optimal coupling of $\mu,\nu$ for $W_2$ takes the form
$$
\Pi=\sum_{i,j}p_{i,j}\delta_{x_i}\otimes\delta_{y_j}\,,\ \ \ i.e.\ \ \  \Pi(x,y)=
\sum_{i,j}p_{i,j}\delta(x-x_i)\delta(y-y_j),
$$
so that
$$
C
_c
=W_2(\mu,\nu)^2.
$$
A general review of the bipartite problem is out of the scope of the present paper, and the reader is referred to the seminal work \cite{mp},  the thesis \cite{gs} which contains an extensive bibliography,  and \cite{BrezisCRAS} for a lucid presentation of the mathematical theory pertaining to this problem. Let us describe the simplest case
$M=N=2$.

In the case of equal masses, that is $m_1=m_2=n_1=n_2=\frac12$, the optimal coupling is shown to be diagonal, in the sense that the mass $\frac12$ is transported from the point $x_1$ to the point $y_1$, and likewise for $x_2$ and $y_2$. 
Thus
$$
\Pi_{op}=\frac12\delta_{x_1}\otimes \delta_{y_1}+\frac12\delta_{x_2}\otimes \delta_{y_2},
$$
or equivalently
$$
\Pi_{op}(x,y)=\frac12\delta(x-x_1)\delta(y-y_1)+\frac12\delta(x-x_2)\delta(y-y_2),
$$
and therefore 
$$
C_c=\frac12(x_1-y_1)^2+\frac12(x_2-y_2)^2\,.
$$

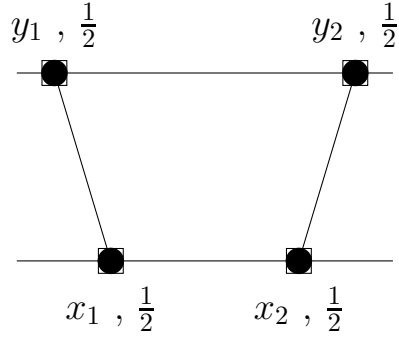
\begin{figure}\label{fig1}
\centering
\begin{tikzpicture}[xscale=2.5,yscale=2.5]
\draw [-] (0,0) -- (2,0);
\draw[-] (0,1) -- (2,1);
\draw[-] (.5,0)--(0.2,1);
\draw[-] (1.5,0)--(1.8,1);
\fill (0.5,0) circle[radius=2pt] node[label=below:$x_1$ ${,}$ $\frac12$,draw]{};
\fill (1.5,0) circle[radius=2pt] node[label=below:$x_2$ ${,}$ $\frac12$,draw]{};
\fill (0.2,1) circle[radius=2pt] node[label=above:$y_1$ ${,}$ $\frac12$,draw]{};
\fill (1.8,1) circle[radius=2pt] node[label=above:$y_2$ ${,}$ $\frac12$,draw]{};
\end{tikzpicture}
\caption{Equal masses} \label{f1}
\end{figure}
\vskip 0.5cm

In the case of unequal masses, let us consider the example where $m_1=\frac{1-\eta}2$ and $m_2=\frac{1+\eta}2$ for some $0<\eta<1$, while $n_1=n_2=\frac12$. In this case, one shows that the optimal transport moves the mass $\frac12$ from $x_2$ to $y_2$, moves the remaining amount of the mass at $x_2$, i.e. $\frac\eta2$, from $x_2$ to $y_1$, and finally moves the mass $\frac{1-\eta}2$ from $x_1$ and $y_1$. Therefore, the optimal coupling in this case is
$$
\Pi_{op}(x,y)=\frac{1-\eta}2\delta(x-x_1)\delta(y-y_1)+\frac\eta2\delta(x-x_2)\delta(y-y_1)+\frac{1}2\delta(x-x_2)\delta(y-y_2),
$$
so that 
$$
C_c
=
\frac{1-\eta}2(x_1-y_1)^2+
\frac\eta2(x_2-y_1)^2+
\frac{1}2(x_2-y_2)^2.
$$

\begin{figure}\label{fig2}
\centering
\begin{tikzpicture}[xscale=2.5,yscale=2.5]
\draw [-] (0,0) -- (2,0);
\draw[-] (0,1) -- (2,1);
\draw[-] (.5,0)--(0.2,1);
\draw[-] (1.5,0)--(1.8,1);
\draw[-] (1.5,0)--(0.2,1);
\fill (0.5,0) circle[radius=2pt] node[label=below:$x_1$ ${,}$ $\frac{1-\eta}2$,draw]{};
\fill (1.5,0) circle[radius=2pt] node[label=below:$x_2$ ${,}$ $\frac{1+\eta}2$,draw]{};
\fill (0.2,1) circle[radius=2pt] node[label=above:$y_1$ ${,}$ $\frac12$,draw]{};
\fill (1.8,1) circle[radius=2pt] node[label=above:$y_2$ ${,}$ $\frac12$,draw]{};
\end{tikzpicture}
\caption{Different masses} \label{f2}
\end{figure}
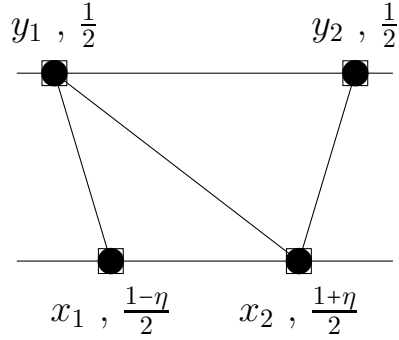

\vskip 1cm 
A quantum analogue to the Wasserstein distance has been recently introduced in \cite{GMouPaul} according to the following general fact.

 When passing from classical to quantum mechanics,

\noindent
1. functions on phase-space should be replaced by operators on the Hilbert space of square integrable functions on the underlying configuration space,  and 

\noindent
2. integration (over phase space) of classical functions should be replaced by the trace of the corresponding operators. Moreover,

\noindent
3. coordinates $q$ of the configuration space should be replaced by the multiplication operator $\hat q$ by the $q$ variable, while momentum coordinates $p$  should be replaced by the operator $\hat p=-i\hbar\nabla$. 

\vskip 1cm
These considerations are  in full accordance with the definition of quantum density matrices as self-adjoint positive operators of trace $1$ on $\fH:=L^2(\bR^d)$. They are also consistent with the definition of couplings $Q$ of two density matrices 
$R$ and $S$ as density matrices on $\fH\otimes\fH$ (identified with $ L^2(\bR^{2d})$) whose  marginals (defined  consistently again as partial traces on the two
factors of $\fH\otimes\fH$) are  equal to $R$ and $S$. In other words
\be\label{quantmarg}\nn
\Tr_{\fH\otimes\fH}((A\otimes I_{\fH})Q)=\Tr_{\fH}(AR),
\quad
\Tr_{\fH\otimes\fH}((I_{\fH}\otimes B)Q)=\Tr_{\fH}(BS)
\ee
for all bounded operators $A,B$ on $\fH$, by analogy with \eqref{clamarg}.

 Moreover they lead naturally to the following definition of the analogue of the Wasserstein distance between two quantum densities $R$ and $S$. Consistently with \eqref{claswass} expressed on the phase-space $\bR^{2d}$, therefore with $n=2d$, we define $MK_2\geq 0$ by
\be\label{quantwass}\nn
MK_2(R,S)^2:=\inf_{Q\text{ coupling of }R\text{ and }S}\Tr{(CQ)},
 \ee
 with
\be\label{couttens}\nn
 C:=(\hat p\otimes \mbox{I}-I\otimes \hat p)^2+
 (\hat q\otimes I-I\otimes \hat q)^2-2d\hbar.
 \ee
 In other words, expressed as an operator on $L^2(\bR^d,dx)\otimes L^2(\bR^d,dy)$,
 \be\label{coutxy}\nn
 C=(x-y)^2-\hbar^2(\nabla_x-\nabla_y)^2-2d\hbar=-4\hbar^2\nabla_{x-y}^2+(x-y)^2-2d\hbar.
 \ee
The operator $\tfrac12(C+2d\hbar)$ is 
a quantum
harmonic oscillator in the variable $(x-y)/\sqrt2$, and this implies in particular that $C=C^*\geq 0$. 

 The quantity $MK_2$ is not a distance as shown in \cite{GMouPaul}  p. 171.
  Nevertheless, it was established in \cite{GMouPaul} the two following links between $MK_2$ and $W_2$\footnote{Note the unessential difference with the definition of the cost $C$ in \cite{GMouPaul,GP1,GP2} created by the shift $-2d\hbar$ and accounts for a shift by $2d\hbar$ in the two next formulas.}. First, for any pair of density matrices $R$ and $S$, the Husimi functions $\widetilde W[R]$ and $\widetilde W[S]$ of $R$ and $S$ satisfy
 \be\label{prophus}\nn
 W_2(\widetilde W[R],\widetilde W[S])^2\leq MK_2(R,S)^2+4d\hbar.
 \ee
On the other hand, if $R$ and $S$ are T\"oplitz operators of symbols $\mu$ and $\nu$,
  \be\label{proptop}
 MK_2(R,S)^2\leq W_2(\mu,\nu)^2.
 \ee
\newcommand{\crochet}[2]{\langle #1| #2\rangle}
Let us recall that a  T\"oplitz operator $T$ (or positive quantization, or anti-Wick ordering quantization)  of symbol a  probability measure $\tau$ 
is\footnote{Here also, we use a different normalization than the one
in \cite{GMouPaul,GP1,GP2}, since we deal exclusively with density matrices. With the present normalization, one has $\Tr{T}=\int_{\bR^{2d}}\tau(dq,dp)$.}
$$
T:=\int_{\bR^{2d}}\ket{q,p}\bra{q,p}\tau(dq,dp),
$$
 where $\ket{q,p}$ is a 
 coherent state at point $(q,p)$ i.e. 
 \begin{equation}\label{CohState}\nn
 \crochet{x}{q,p}:=
 (\pi\hbar)^{-d/4}e^{-(x-q)^2/2\hbar}e^{ipx/\hbar}.
\end{equation}
We also recall the definition of the Husimi function of a density matrix $R$:
  $$
  \widetilde W[R](q,p):=(2\pi\hbar)^{-d}\bra{q,p}R\ket{q,p}.
  $$
\vskip 0.5cm
The functional $MK_2^2$ (more precisely $MK_2^2+2d\hbar$ with the definition chosen in the present paper) has been systematically used and extended in \cite{GMouPaul,GP1,GP2} in order to study various problems, such as the validity of
the mean-field limit uniformly in $\hbar$, the semiclassical approximation of quantum dynamics, and the problem of metrizing of the set of quantum densities in the semiclassical regime.

\smallskip
Given the importance of optimal transport in the field of statistics in the problem of comparing probability measures, there have been various attempts at defining analogues of the Wasserstein, or Monge-Kantorovich distances in the quantum 
setting. For instance, the reference \cite{ZyckSlom} proposed to consider the original Monge distance (also called the Kantorovich-Rubinstein distance, or the Wasserstein distance of exponent $1$) between the Husimi transforms of the 
density operator. However, propagating this distance with the usual quantum dynamics may not be easy, because the dynamics of the Husimi transform of a density operator by the von Neumann equation is quite involved \cite{AthaMausPaul}. 
A big advantage of the quantity $MK_2$ introduced in \cite{GMouPaul} is that it is directly defined in terms of density operators, and therefore easily propagated by the usual quantum dynamics, including $N$-body quantum dynamics, for 
which it has been defined originally. Besides the quantity $MK_2$ appeared in \cite{GMouPaul}, other quantum analogues of the Wasserstein distance of exponent $2$ have been proposed by several other authors. For instance a quantum 
analogue of the Benamou-Brenier formula (see Theorem 8.1 in chapter 8 of \cite{VillaniTOT}) for the classical Wasserstein distance of exponent $2$ is studied in detail in \cite{CarlMaas1,CarlMaas2}, and this idea has been used to obtain 
a quantum equivalent of the so-called HWI inequality \cite{RouzDatt}. More recently, other propositions for generalizing Wasserstein distances to the quantum setting have emerged, such as \cite{Ikeda} (which seems essentially focussed 
on pure states) or \cite{PalmTrev}, which is very close to our definition of $MK_2$, except that the set of couplings used in the minimization is different.

\vskip 0.5cm
The quantum bipartite problem can be therefore stated as follows, in close analogy with the classical picture introduced earlier.

One considers two density matrices built in terms of the positions and masses already used for the classical bipartite problem, in the following way
$$
R
=
\sum_{i=1}^Mm_i\ket{x_i,0}\bra{x_i,0},\qquad S=\sum_{j=1}^Nn_j\ket{y_j,0}\bra{y_j,0}.
$$

Indeed, it is natural to associate coherent states to material points, as they saturate the Heisenberg uncertainty inequalities. Moreover, one sees that $R$ and $S$ are precisely the T\"oplitz operators of symbols $\mu$ and $\nu$ respectively.

The quantum bipartite problem consists then in finding an optimal coupling of $R$ and $S$ for $MK_2(R,S)$ and the optimal quantum cost defined as
$$
C_q:=
MK_2(R,S).
$$
Since $R$ and $S$ are T\" oplitz operators, we know from \eqref{proptop} that
$$
C_q\leq C_c.
$$
The question we address in this paper is whether there exist pairs of density matrices for which
$$
C_q< C_c.
$$
In other words, we address the question of whether quantum optimal transportation can be cheaper than its classical analogue.
\vskip 1cm

We shall study the two cases
introduced at the beginning of this section and described in Figures 1 and 2. For the sake of simplicity, we shall take $x_1=-x_2=-a, y_1=-y_2=-b$, with $a<b$ in the equal mass case, and $a=b$ in the unequal mass case.

In the equal mass case, studied in  Section \ref{eqmass}, both classical and quantum transport are achieved without splitting mass for each particle: the two costs are shown to be equal (see \eqref{eqm}), and an optimal quantum coupling is 
the T\" oplitz quantization to the optimal classical coupling.

In Section \ref{diffmass} we study the case 
of different masses and construct a family of examples for which the optimal quantum cost is strictly cheaper than the classical one (see \eqref{diffm}). 

In addition, we show in Section \ref{conc} that an optimal quantum coupling is not always the T\"oplitz quantization of a classical coupling. This can be rephrased by saying that  an optimal quantum transport can be different from the natural quantization of any underlying classical transport. In fact, in the unequal mass case treated in this paper, no quantum optimal transport  corresponds to a classical transport, optimal or not: they all involve strictly quantum effects.




\section{The equal mass case}\label{eqmass}
For $a,b>0$ we will  transport a superposition of two density matrices which are pure states associated  to two coherent states of null momenta localized at $+a$ and $-a$ towards a similar density matrix associated to the points $(\pm b,0)$
in phase space. In other words, we consider the coherent states denoted $\ket{c}$ for simplicity (instead of $\ket{c,0}$,
i.e. $
\crochet{x}{c}:=(\pi\hbar)^{-1/4}e^{-(x-c)^2/2\hbar}
$)
and consider the two density matrices
$$
R:=\frac12(\ket{a}\bra{a}+\ket{-a}\bra{-a}), \qquad S:=\frac12(\ket{b}\bra{b}+\ket{-b}\bra{-b}).
$$
Define
$$
\lambda:=\bra{a}{-a}\rangle=e^{-a^2/\hbar},\qquad
\mu:=\bra{b}{-b}\rangle=e^{-b^2/\hbar},
$$
and consider the two pairs of orthogonal vectors
\be\label{phipsi}
\phi_\pm:=\frac{\ket{a}\pm\ket{-a}}{\sqrt{2(1\pm\lambda)}},\qquad
\psi_\pm:=\frac{\ket{b}\pm\ket{-b}}{\sqrt{2(1\pm\mu)}}.
\ee
Hence
$$
R=\alpha_+\ketbra{\phi_+}{\phi_+}
+
\alpha_-\ketbra{\phi_-}{\phi_-},\qquad
S=\beta_+\ketbra{\psi_+}{\psi_+}
+
\beta_-\ketbra{\psi_-}{\psi_-},
$$
with
$$
\alpha_+:=\frac12(1+\lambda),\qquad \alpha_-=\frac12(1-\lambda),\qquad
\beta_+:=\frac12(1+\mu),\qquad\beta_-=\frac12(1-\mu).
$$

In  the whole present paper, we will only use  couplings of $R$ and $S$ that act from the four-dimensional linear span of $\phi_\pm\otimes \psi_\pm$ to itself. Therefore, in order to compute $\Tr{(CQ)}$ for such couplings, we need to project the cost operator
$C$ on the basis $\{\phi_+\otimes\psi_+,\phi_+\otimes \psi_-,\phi_-\otimes\psi_+,\phi_-\otimes\psi_-\}$. This is a tedious but straightforward computation which results in the following $4\times 4$ matrix:
\be\label{cost}
C=\begin{pmatrix}
\cA&0&0&\gamma\\
0&\cB&\delta&0\\
0&\delta&\cC&0\\
\gamma&0&0&\cD
\end{pmatrix}.
\ee
where 
$$
\begin{aligned}
\cA=a^2\frac{1-\lambda}{1+\lambda}+b^2\frac{1-\mu}{1+\mu}&,\qquad\cB=a^2\frac{1-\lambda}{1+\lambda}+b^2\frac{1+\mu}{1-\mu}, \qquad
\gamma=-\frac{2ab(1-\lambda\mu)}{\sqrt{(1-\lambda^2)(1-\mu^2)}},
\\
\cC=a^2\frac{1+\lambda}{1-\lambda}+b^2\frac{1-\mu}{1+\mu}&,\qquad\cD=a^2\frac{1+\lambda}{1-\lambda}+b^2\frac{1+\mu}{1-\mu}, \qquad\delta=-\frac{2ab(1+\lambda\mu)}{\sqrt{(1-\lambda^2)(1-\mu^2)}}.
\end{aligned}
$$


As a warm up in order to find an ansatz for the general case, let us  first neglect the contributions of $\lambda,\mu$, exponentially small in the Planck constant. In this case $\alpha_\pm=\beta_\pm=\frac12$, and the cost is equal to
$$
C_0=\begin{pmatrix}
a^2+b^2&0&0&-2ab\\
0&a^2+b^2&-2ab&0\\
0&-2ab&a^2+b^2&0\\
-2ab&0&0&a^2+b^2
\end{pmatrix},
$$
On the other hand, one has
$$
Q_0:=\frac14\begin{pmatrix}
1&0&0&1\\
0&1&1&0\\
0&1&1&0\\
1&0&0&1
\end{pmatrix}\geq 0,
$$
since the spectrum of $Q_0$ is easily shown to be $\{0,\tfrac12\}$ by using the elementary formula
\be\label{det}
\det\begin{pmatrix}
\bar a&0&0&\gamma\\
0&\bar b&\delta&0\\
0&\delta&\bar c&0\\
\gamma&0&0&\bar d
\end{pmatrix}=(\bar a\bar d-\gamma^2)(\bar b\bar c-\delta^2)\mbox{ for all } \bar a,\bar b,\bar c,\bar d,\gamma,\delta. 
\ee
Moreover, one easily checks that $\Tr_2{Q_0}=R$ and $\Tr_1{Q_0}=S$ so that $Q_0$ is  a coupling of $R$ and $S$.

Another easy computation shows that
$$
\Tr{(CQ_0)}=(a-b)^2\,.
$$

Therefore
 $$
MK_2(R,S)^2\leq (a-b)^2=W_2(\tfrac12(\delta_{-a}+\delta_a),\tfrac12(\delta_{-b}+\delta_b))^2.
$$
\vskip 0.5cm

For the ``true" case $\lambda,\mu\not=0$, we make the following ansatz on the coupling $Q$
$$
Q=Q_0+\frac14
\begin{pmatrix}
p+\lambda+\mu&0&0&u\\
0&-p+\lambda-\mu&v&0\\
0&v&-p-\lambda+\mu&0\\
u&0&0&p-\lambda-\mu
\end{pmatrix},\  p,u,v\in\bR.
$$
Straightforward computations show that
$$
\Tr Q=\Tr Q_0=1,\qquad \Tr_2Q=\Tr_2Q_0=R,\qquad \Tr_1Q=\Tr_1Q_0=S.
$$
Using again \eqref{det} shows that 
$$
Q\geq 0 \Longleftrightarrow -1+\sqrt{(\lambda+\mu)^2+(1+u)^2}
\leq p
\leq
1-\sqrt{(\lambda-\mu)^2+(1+v)^2}.
$$
Therefore, assuming that $p,u,v$ satisfy this constaint, $Q$ is a coupling of $R$ and $S$.

\vskip 1cm
Denoting $U:=1+u$ and $V:=1+v$, we compute $W:=\Tr(CQ)$ by using \eqref{cost}: 
\begin{eqnarray}
4W
&=&
2\gamma U+2\delta V+ p(\cA-\cB-\cC+\cD)+\cA+\cB+\cC+\cD\nn
\\
&&+(\lambda+\mu)(\cA-\cD)+(\lambda-\mu)(\cB-\cC)\nn
\\
&=&
2\gamma U+2\delta V+ p(\cA-\cB-\cC+\cD)+W'=2\gamma U+2\delta V+W',
\nn
\end{eqnarray}
with
\begin{eqnarray}
W'&:=&\cA+\cB+\cC+\cD+\lambda(
\cA+\cB-\cC-\cD
)+\mu(
\cA-\cB+\cC-\cD
)
\nn\\
&=&
4\left(a^2\frac{1+\lambda^2}{1-\lambda^2}+b^2\frac{1+\mu^2}{1-\mu^2}\right)
-8a^2\frac{\lambda^2}{1-\lambda^2}
-8b^2\frac{\mu^2}{1-\mu^2}\nn
\\
&=&
4(a^2
+b^2\nn
).
\end{eqnarray}
Since $W$ is linear in $U,V$, we minimize $\gamma U+\delta V$ by taking
$$
U=\sqrt{(p+1)^2-(\lambda+\mu)^2}\qquad\text{ and }\qquad
V= \sqrt{(p-1)^2-(\lambda-\mu)^2},
$$
and, since $\delta\leq\gamma$, we conclude that 
$$
4W=2T+W',
$$ 
where
$$
\begin{aligned}
T=&-\max_{-1+\lambda-\mu\leq p\leq 1-(\lambda-\mu)}\left(-\gamma\sqrt{(p+1)^2-(\lambda+\mu)^2}-\delta\sqrt{(p-1)^2-(\lambda-\mu)^2}\right)
\\
=&\frac{-2ab}{\sqrt{(1-\lambda^2)(1-\mu^2)}}\max_{-1+\lambda-\mu\leq p\leq 1-(\lambda-\mu)}\left((1-\lambda\mu)\sqrt{(p+1)^2-(\lambda+\mu)^2}\right.
\\
&\qquad\qquad\qquad\qquad\qquad\qquad\qquad\qquad\left.+(1+\lambda\mu)\sqrt{(p-1)^2-(\lambda-\mu)^2}\right).
\end{aligned}
$$
One can check that the max is attained for $p=\lambda\mu\to 0$ as $\hbar\to 0$, and that
$$
T=
-\frac{4ab}{\sqrt{(1-\lambda^2)(1-\mu^2)}}
\sqrt{1+\lambda^2\mu^2-\lambda^2-\mu^2}=-4ab.
$$ 
Eventually, we arrive at the 
same result as in the semiclassical regime $\lambda=\mu=0$, viz.
\be\label{estinf}
MK_2(R,S)^2
\leq
(a-b)^2.
\ee
Since $R$ and $S$ are T\"oplitz operator, the inequality \eqref{estinf} was already known by using \eqref{proptop}. Nevertheless we gave this explicit computation as we believe the result to be valid for more general density matrices.

\vskip 1cm
In order to get a lower bound for $MK_2(R,S)$, we shall use a dual version of the definition of $MK_2$, proved in \cite{cgp},
that is a quantum version of the Kantorovitch 
duality theorem for $W_2$ (see \cite{VillaniTOT,VillaniTOT2}):
\be\label{kant}\nn
MK_2(R,S)^2=\sup_{
A=A^*,\ B=B^*\text{ bounded operators on }\fH\atop \text{ such that }A\otimes I+I\otimes B\le C}
\Tr(RA+SB).
\ee
\vskip 1cm
We make the following diagonal ansatz on $A$ and $B$:
$$
A
=
\begin{pmatrix}
\alpha_1&0\\0&\alpha_2
\end{pmatrix}
\qquad
B=
\begin{pmatrix}
\beta_1&0\\0&\beta_2
\end{pmatrix}\,,
$$
so that 
$$
A\otimes I
=
\begin{pmatrix}
\alpha_1&0&0&0\\
0&\alpha_1&0&0\\
0&0&\alpha_2&0\\
0&0&0&\alpha_2
\end{pmatrix}
\quad\text{ 
and
 }\quad
I\otimes B
=
\begin{pmatrix}
\beta_1&0&0&0\\0&\beta_2&0&0\\0&0&\beta_1&0\\
0&0&0&\beta_2
\end{pmatrix}\,.
$$
Hence
$$
A\otimes I+I\otimes B-C:=\begin{pmatrix}
\bar a&0&0&-\gamma\\
0&\bar b&-\delta&0\\0&-\delta&\bar c&0\\-\gamma&0&0&\bar d
\end{pmatrix}\,,
$$ 
and, according to \eqref{cost},
$$
\bar a=\alpha_1+\beta_1-\cA,\quad\bar b=\alpha_1+\beta_2-\cB,\quad\bar c=\alpha_2+\beta_1-\cC,\quad\bar d=\alpha_2+\beta_2-\cD.
$$
Notice that
$$
\bar a+\bar d=\bar b+\bar c.
$$
Using \eqref{det} to compute the characteristic polynomial of $A\otimes I+I\otimes B-C$, we find that
\be\label{cons}
A\otimes I+I\otimes B\leq C\Longleftrightarrow 
\bar a+\bar d\leq-\sqrt{(\bar a-\bar d)^2+4\gamma^2}
\mbox{ and }\bar b+\bar c\leq-\sqrt{(\bar b-\bar c)^2+4\delta^2}.
\ee

Moreover,
\begin{eqnarray}
\Tr(AR+BS)&=&
%
\frac12(\alpha_1+\alpha_2+\beta_1+\beta_2)
+
\frac\lambda2(\alpha_1-\alpha_2)
+
\frac\mu2(\beta_1-\beta_2)\nn\\
&
=&
\frac14(\bar a+\bar b+\bar c+\bar d)
+
\frac14(\bar a+\bar b-\bar c-\bar d)\lambda\\
&+&
\frac14(\bar a-\bar b+\bar c-\bar d)\mu
+
a^2+b^2\nn.
\end{eqnarray}

Let us denote 
$$
x:=\bar a+\bar d=\bar b+\bar c,
$$
so that
\be\label{traceab}
\Tr(AR+BS)=\frac12x+
\frac14(\lambda+\mu)(\bar a-\bar d)
+
\frac14(\lambda-\mu)(\bar b-\bar c)+a^2+b^2.
\ee

The constraints \eqref{cons} are expressed as
$$
\begin{aligned}
x=&\,\bar a+\bar d&&\leq-\sqrt{(\bar a-\bar d)^2+4\gamma^2}\,,
\\
x=&\,\bar b+\bar c&&\leq-\sqrt{(\bar b-\bar c)^2+4\delta^2}\,.
\end{aligned}
$$
Without loss of generality we assume that $\lambda\geq\mu$, that is to say $a<b$. Since the right hand side of \eqref{traceab} is linear in $x$, in $(\bar a-\bar d)$, and in $(\bar b-\bar c)$, one has to saturate the constraints to maximize 
$\Tr(AR+BS)$. In other words, we must take
$$
\bar a-\bar d=\sqrt{x^2-4\gamma^2},\quad\text{ and }\quad  \bar b-\bar c=\sqrt{x^2-4\delta^2}.
$$
Since $\delta\leq\gamma\leq0$, this amounts to computing
$$
\max_{x\leq2\delta}f(x),\qquad\text{ with }f(x):=\frac x2+\frac14(\lambda+\mu)\sqrt{x^2-4\gamma^2}+\frac14(\lambda-\mu)\sqrt{x^2-4\delta^2}\,.
$$
We check that $f'(x)$ is an increasing function of $x^2$, so that the maximum of $f(x)$ for $x\le 2\delta$ is attained at
$$
f'(x)=0\iff x=-\frac{4ab(1-\lambda^2\mu^2)}{(1-\lambda^2)(1-\mu^2)},\qquad\text{ which implies }f(x)=-2ab.
$$
We conclude from \eqref{traceab} that
$$
MK_2(R,S)^2\geq\Tr{(AR+BS)}\geq a^2+b^2-2ab=(a-b)^2\,.
$$
Together with \eqref{estinf}, this implies that
$$
MK_2(R,S)^2=(a-b)^2=W_2(\tfrac12(\delta_{-a}+\delta_a),\tfrac12(\delta_{-b}+\delta_b))^2.
$$
Therefore,
\be\label{eqm}
C_q=C_c,
\ee
so that the classical and the quantum optimal transport costs are equal in this case.
\vskip 1cm

\section{The unequal mass case}\label{diffmass}
In this section, we construct a family of density matrices $R$ and $S$ for which the quantum cost of optimal transport is smaller than the classical analogous cost.

With the same notations as in previous section, we set
$$
R:=\tfrac{1+\eta}2\ket{a}\bra{a}+\tfrac{1-\eta}2\ket{-a}\bra{-a}, \qquad S:=\tfrac12\ket{a}\bra{a}+\tfrac12\ket{-a}\bra{-a},\qquad 0<\eta<1.
$$
In other words, we consider the same situation as in the previous section with $a=b$, but with different masses for the quantum density matrix $R$.

In the orthonormal basis $\{\phi_+,\phi_-\}$, the density matrix $R$ takes the form
\be\nn
R=\begin{pmatrix}
\frac{1+\lambda}2&\frac\eta 2\sqrt{1-\lambda^2}\\	\\
\frac\eta 2\sqrt{1-\lambda^2}\quad&\frac{1-\lambda}2
\end{pmatrix}\,,
\ee
while $S$ is the same as before.

We define the ``quantized classical'' coupling as
\be \label{notpointvirgule}
Q_c:=\tfrac1 2\ket{a;a}\bra{a;a}
+
\tfrac{1-\eta}2\ket{-a;-a}\bra{-a;-a}
+
\tfrac\eta 2\ket{a;-a}\bra{a;-a},
\ee
with the obvious notation
$$
\bra{a;b}:=\bra{a}\otimes\bra{b};\qquad\ket{a;b}:=\ket{a}\otimes\ket{b}.
$$

Obviously $Q_c\geq 0$ by construction, and 
$$
\Tr_2(Q_c)=\tfrac12\ket{a}\bra{a}+\tfrac{\eta}2\ket{a}\bra{a}+\tfrac{1-\eta}2\ket{-a}\bra{-a}=R,\quad\text{ while }\Tr_1(Q_c)=S.
$$

Viewed as a matrix in the basis $\{\phi_+\otimes\psi_+,\phi_+\otimes \psi_-,\phi_-\otimes\psi_+,\phi_-\otimes\psi_-\}$, 
\be\label{qcmat}
Q_c=
\small{
\begin{pmatrix}
\frac14(1+\lambda)^2&0&\frac14\eta\sqrt{1-\lambda}(1+\lambda)^{\frac32}&\frac14(-1+\eta)(-1+\lambda^2)\\
0&\frac14(1-\lambda^2)&\frac14(-1+\eta)(-1+\lambda^2)&\frac14\eta(1-\lambda)^{\frac32}\sqrt{1+\lambda}\\
\frac14\eta\sqrt{1-\lambda}(1+\lambda)^{\frac32}&\frac14(-1+\eta)(-1+\lambda^2)&\frac14(1-\lambda^2)&0\\
\frac14(-1+\eta)(-1+\lambda^2)&\frac14\eta(1-\lambda)^{\frac32}\sqrt{1+\lambda}&0&\frac14(-1+\lambda)^2
\end{pmatrix}.
}
\ee
With \eqref{cost}, we easily compute 
%
\be\label{qclass}
\Tr{(CQ_c)}=2\eta a^2=W_2(\tfrac{1+\eta}2\delta_{a}+\tfrac{1-\eta}2\delta_{-a}, \tfrac12\delta_a+\tfrac12\delta_{-a})^2.
\ee
Indeed, let us recall the classical optimal transport from $R$ to $S$ in this case: first, one ``moves'' the amount of mass $\frac12$ from $a$ in $R$ to $a$ in $S$. The amount of mass $\tfrac{\eta}2$ remaining at $a$ in $R$ is transported to $-a$ 
in $S$, and the outstanding amount of mass $\frac{1-\eta}2$, located at $-a$ in $R$, is ``transported'' to $-a$ in $S$ (see Figure {2}).
 

 For each $\epsilon>0$, set
 \be \label{qeps}
 Q_\epsilon:=Q_c+\epsilon Q_q,
 \ee
 with 
 \be\nn\
   Q_q:=\begin{pmatrix}
 1&0&0&-1\\
 0&-1&1&0\\
 0&1&-1&0\\-1&0&0&1
 \end{pmatrix}
\ee
One easily checks that
\be\nn 
\Tr_1{(Q_q)}=\Tr_2{(Q_q)}=\Tr{(Q_q)}=0,
\ee
so that
\be\label{traceqq}
\Tr_1{(Q_\epsilon)}=S,\quad\text{ and }\quad\Tr_2{(Q_\epsilon)}=R,\quad\text{ so that }\Tr{(Q_\epsilon)}=1.
\ee

The characteristic polynomial of $Q_c$ is found to be of the form
$$
\det{(Q_c-tI)}=tP_3(t),
$$
where $P_3$ is a cubic polynomial satisfying 
$$
P_3(0)=-\tfrac\eta 8(1-\eta)(1-\eta^2)<0.
$$
Therefore the spectrum of $Q_c$ is $\{0,\lambda_1>0,\lambda_2>0,\lambda_3>0\}$ since $Q_c=Q^*_c\ge 0$. One can also check that
\be \label{mnb1}
\det{(Q_\epsilon-tI)}|_{t=0}=\det{Q_\epsilon}=\epsilon\eta\lambda^2(1-\eta)(1-\lambda^4)+O(\epsilon^2)>0\quad\text{ for }0<\epsilon\ll 1,
\ee
together with
$$
\frac d{dt} \det{(Q_c-tI)}|_{t=0}:=P_3(0)<0.
$$
Hence there exists $D$ (independent of $\epsilon$) such that
\be \label{mnb2}
\frac d{dt}\det{(Q_\epsilon-tI)}|_{t=0}\le D<0\qquad\text{ for }0<\epsilon\ll 1.
\ee
Both \eqref{mnb1} and \eqref{mnb2} clearly imply that $\det{(Q_\epsilon-tI)}$ has a positive zero that is $\epsilon$-close to $0$, and three other roots which are $\epsilon$-close to $\lambda_1$, $\lambda_2$ and $\lambda_3>0$ respectively.
Therefore, $Q_\epsilon=Q_\epsilon^*>0$ for $0<\epsilon\ll 1$, and \eqref{traceqq} implies that $Q_\epsilon$ is a coupling of $R$ and $S$.
 
Another elementary computation shows that
\be\nn
\Tr{(CQ_q)}=-\frac{8a^2\lambda^2}{1-\lambda^2},
\ee
so that
$$
\begin{aligned}
MK_2(R,S)^2\le&\Tr{(CQ_\epsilon)}=\Tr{(CQ_c)}-\epsilon\frac{8a^2\lambda^2}{1-\lambda^2}
\\	\\
<& W_2(\tfrac{1+\eta}2\delta_{a}+\tfrac{1-\eta}2\delta_{-a}, \tfrac12\delta_a+\tfrac12\delta_{-a})^2,
\end{aligned}
$$
for each $\epsilon$ satisfying $0<\epsilon\ll 1$, according to formula \eqref{qclass}. 
In other words, 
\be\label{diffm}
C_q<C_c,
\ee
the quantum cost 
is (strictly) below the classical cost.

\section{Concluding remarks on quantum optimal transport}\label{conc}

The result of Section \ref{eqmass} shows that, in the equal mass case, an optimal coupling is given by the following matrix in the basis $\{\phi_+\otimes\psi_+,\phi_+\otimes \psi_-,\phi_-\otimes\psi_+,\phi_-\otimes\psi_-\}$:
$$
Q
=
\tiny{
\frac14
\begin{pmatrix}
1+\lambda\mu+\lambda+\mu&0&0&\sqrt{(1+\lambda\mu)^2-(\lambda+\mu)^2}\\
0&1-\lambda\mu+\lambda-\mu&\sqrt{(1-\lambda\mu)^2-(\lambda-\mu)^2}&0\\
0&\sqrt{(1-\lambda\mu)^2-(\lambda-\mu)^2}&1-\lambda\mu-\lambda+\mu&0\\
\sqrt{(1+\lambda\mu)^2-(\lambda+\mu)^2}&0&0&1+\lambda\mu-\lambda-\mu
\end{pmatrix}.
}
$$
Using  \eqref{phipsi} and with the same notation as in \eqref{notpointvirgule}, the optimal coupling $Q$ can be put in the form
\be\label{qab}
 Q=\frac12\left(\ket{a;b}\bra{a;b}
+
\ket{-a;-b}\bra{-a;-b}
\right).
\ee
In other words, $Q$ is the T\"oplitz operator  of symbol
$$
\Pi(q,p;q',p')=
\tfrac12
\delta_{(-a,0)}(q,p)\delta_{(-b,0)}(q',p')+\frac12
\delta_{(a,0)}(q,p)\delta_{(b,0)}(q',p').
$$
Likewise, we recall that $R$ is the T\"oplitz operator of symbol
\be\label{topr}\nn
\mu(q,p)=\tfrac12\big(\delta_{(-a,0)}(q,p)+\delta_{(a,0)}(q,p)\big),
\ee
while 
$S$ is the T\"oplitz operator of symbol
\be\label{tops}\nn
\nu(q,p)=\tfrac12\big(\delta_{(-b,0)}(q,p)+\delta_{(b,0)}(q,p)\big).
\ee


Therefore, 
\begin{eqnarray}\label{ksbtop}
\Pi(q,p;q',p')
&=&
\tfrac12\big(
\big(\delta_{(-a,0)}(q,p)+\delta_{(a,0)}(q,p)\big)
\delta((q',p')-\Phi(q,p))
\big)\nn\\
&=&
\mu(q,p)\delta((q',p')-\Phi(q,p)),
\end{eqnarray}
where $\Phi$ is any map satisfying $\Phi(a,0)=(b,0)$ and $\Phi(-a,0)=(-b,0)$.

The second equality in \eqref{ksbtop} 
says the following: in the equal mass case, in agreement with  \eqref{ksb}, an optimal quantum coupling $Q$ is the T\"oplitz operator of symbol the classical optimal coupling associated to the optimal transport map
$$
\big((-a,0),(a,0)\big)\mapsto\big((-b,0),(b,0)\big).
$$
\vskip 0.5cm

In the unequal mass case treated in Section \ref{diffmass}, the coupling $Q_c$ defined by \eqref{notpointvirgule} is also a T\"oplitz operator, with symbol
\begin{eqnarray}\label{qctop}\nn
\Pi_c(q,p;q',p')&=&
\tfrac12\delta_{(a,0)}(q,p)\delta_{(a,0)}(q',p')\nn\\
&+&
\tfrac{1-\eta}2\delta_{(-a,0)}(q,p)\delta_{(-a,0)}(q',p')
+
\tfrac\eta2\delta_{(a,0)}(q,p)\delta_{(-a,0)}(q',p').\nn
\end{eqnarray}

This expression is easily interpreted as the optimal coupling associated to the ``transport'' introduced in Section \ref{intro}, Figure 2, exactly as in the equal mass case. But, as explained in the previous section, $Q_c$ cannot be an optimal 
coupling, since the coupling $Q_\epsilon$ defined  by \eqref{qeps} leads to a strictly lower quantum cost. 

We did not compute any optimal coupling in this situation. Observe however that, thanks to \eqref{qcmat} and \eqref{phipsi} specialized to $a=b$ (so that $\lambda=\mu$), one can expand  $Q_q$ in  the form
$$
Q_q=\sum_{i,j,k,l=\pm1}q_{i,j,k,l}\ket{ia;ja}\bra{ka;la}.
$$
The contribution of the ``diagonal" terms $q_{i,j,i,j}$ defines a T\"oplitz operator, unlike the off-diagonal terms such as $q_{1,1,-1,1}=\frac{-4\lambda}{(1-\lambda^2)^2}\neq 0$ for instance.

In general, when $R$ and $S$ are T\"oplitz operators of symbols $\mu$ and $\nu$ satisfying $MK_2(R,S)<W_2(\mu,\nu)$, no optimal coupling $Q_{op}$ of $R$ and $S$ can be a T\"oplitz operator: if such was the case, the T\"oplitz 
symbol of $Q_{op}$ would be a coupling of $\mu$ and $\nu$ with classical transport cost $MK_2(R,S)<W_2(\mu,\nu)$, which is impossible. The presence of nonclassical off-diagonal terms in $Q_{op}$, such as 
$q_{1,1,-1,1}=\frac{-4\lambda}{(1-\lambda^2)^2}\neq 0$ in the example discussed above, are precisely the reason why quantum optimal transport can be cheaper in this case than classical optimal transport.

Finally, observe that both $W_2(\tfrac{1+\eta}2\delta_{a}+\tfrac{1-\eta}2\delta_{-a}, \tfrac12\delta_a+\tfrac12\delta_{-a})^2-\Tr{(CQ_\epsilon)}$ and $q_{1,1,-1,1}$ are exponentially small as $\hbar\to 0$, but of course are not small for $\hbar=1$.

\vskip 1cm
{\bf Acknowledgements}.  This work has been partially carried out thanks to the support of the LIA AMU-CNRS-ECM-INdAM Laboratoire Ypatie des Sciences Math\'ematiques (LYSM).
\vskip 1cm
%
%
%
%
%
%
%
%
%
%

\end{document}